\RequirePackage{ifpdf}
\ifpdf 
\documentclass[pdftex]{sigma}
\else
\documentclass{sigma}
\fi

\begin{document}


\allowdisplaybreaks

\renewcommand{\PaperNumber}{069}

\FirstPageHeading

\renewcommand{\thefootnote}{$\star$}

\ShortArticleName{Yangian of the Strange Lie Superalgebra of $Q_{n-1}$ Type, Drinfel'd Approach}

\ArticleName{Yangian of the Strange Lie Superalgebra of $\boldsymbol{Q_{n-1}}$ Type,
Drinfel'd Approach\footnote{This paper is a
contribution to the Vadim Kuznetsov Memorial Issue `Integrable
Systems and Related Topics'. The full collection is available at
\href{http://www.emis.de/journals/SIGMA/kuznetsov.html}{http://www.emis.de/journals/SIGMA/kuznetsov.html}}}

\Author{Vladimir STUKOPIN}

\AuthorNameForHeading{V. Stukopin}

\Address{Don State Technical University, 1 Gagarin Square, Rostov-na-Donu, 344010 Russia}
\Email{\href{mailto:stukopin@math.rsu.ru}{stukopin@math.rsu.ru}}

\ArticleDates{Received November 01, 2006, in f\/inal form May
06, 2007; Published online May 22, 2007}

\Abstract{The Yangian of the strange Lie superalgebras in Drinfel'd
realization is def\/ined. The current system generators and def\/ining relations are described.}

\Keywords{Yangian; strange Lie superalgebra;  Drinfel'd
realization; Hopf structure; twisted current bisuperalgebra}

\Classification{17B37}

\section{Introduction}
In this paper, following Drinfel'd, we def\/ine
the Yangian of the strange Lie superalgebra of  $Q_{n-1}$ type.
Recall that the Yangian of the simple Lie algebra was def\/ined by V.~Drinfel'd
as a~quantization of the polynomial currents Lie bialgebra (with values
in this simple Lie algebra) with coalgebra structure def\/ined by rational Yang $r$-matrix
\cite{Dr,Dr1,Dr2,Ch-Pr}. The Yangian of the reductive Lie algebra
can be given the same def\/inition in special cases. The object dual to Yangian
(of the general linear Lie algebra $gl(n)$))  was studied before by L.~Faddeev
and others while working  on Quantum Inverse Scatering Method (QISM).
We call this def\/inition of Yangian the RFT approach. V.~Drinfel'd shows this object to be isomorphic to
the Yangin of $\mathfrak{gl}(n)$. It is the RFT approach
that is usual for papers devoted to Yangians and def\/ines the Yangian
as the algebra generated by matrix elements of Yangians irreducible
representations according to Drinfel'd (see \cite{Mol} and \cite{Zh1,Zh2,Zh3}
for Yangians of the Lie superalgebras). More precisely, the Yangian
can be viewed as the Hopf algebra generated by matrix elements of a matrix $T(u)$
(so-called transfer matrix) with the commutation def\/ining relations:
\[
R(u-v)T_1(u)T_2(v)= T_2(v)T_1(u)R(u-v),
\]
where $T_1(u)=T(u)\otimes E$, $T_2(v)=E\otimes T(v)$, $E$ is an identity matrix $R(u)$
is some rational matrix-function (with the values in ${\rm End}\,(V\otimes V$)).
$R(u)$ is called a quantum $R$-matrix (see \cite{Mol}).

The RFT approach is used in the paper \cite{N} (see also \cite{Zh1,Zh2,Zh3})
to def\/ine the Yangian of basic Lie superalgebra, while  the Drinfel'd's
one is used in papers \cite{St,St1,St3}. We cannot use directly Drinfel'd
approach for def\/ining the Yangian of strange Lie superalgebra~\cite{K}
as this Lie superalgebra does not have nonzero invariant bilinear forms.
Hence we cannot def\/ine the Lie bisuperalgebra's
structure on the polynomial currents Lie superalgebra with values
in the strange Lie superalgebra. M.~Nazarov noted that this structure
can be def\/ined on the twisted current Lie superalgebra and used the RFT approach to do it (see~\cite{N1}).
In this paper we use Nazarov's idea to def\/ine the Yangian
of strange Lie superalgebra according to Drinfel'd.
Our def\/inition can be used for further research of the Yangian
of strange Lie superalgebra. It is quite convenient to study such
problems as exact description of the quantum double of Yangian
of strange Lie superalgebra and computation of Universal $R$-matrix of quantum double.
Our method can also be extended on the other twisted current Lie algebras and Lie
superalgebra quantization description.  As a result, we obtain so-called the twisted
Yangians without a Hopf (super)algebra  structure (in general) but being comodule over some Hopf (super)algebras.

 Note that the Yangian of the basic Lie superalgebra $A(m,n)$ was
 def\/ined according to Drinfel'd in \cite{St} where the  Poincar\'e--Birkgof\/f--Witt theorem
 (PBW-theorem) and the theorem on existence of pseudotriangular
 structure on Yangian (or the theorem on existence of the universal $R$-matrix)
 are also proved. Further, in \cite{St1,St3} the quantum double
 of the Yangian of the Lie superalgebra $A(m,n)$ is described,
 and the multiplicative formulas for the universal  $R$-matrices
 (for both quantum double of the Yangian and Yangian) are obtained.
 This paper is the consequence of \cite{St,St2}
 and extension of some of their results on the Yangian of the ``strange'' Lie superalgebra.

 Following Drinfel'd, we def\/ine the Yangian of the strange Lie
 superalgebra of $Q_{n-1}$ type and describe the current system
 of Yangian generators and def\/ining relations, which is an analogue
 of the same  system from \cite{Dr2}.   The problem of the equivalency between our
 def\/inition and Nazarov's one is not resolved yet.
 The problem of constructing of the  explicit formulas
 def\/ining the isomorphism between the above realization (as in \cite{Cr} in the case of $sl_n$)
 is open and seems very interesting. This will be discussed in further papers.

\section{Twisted current bisuperalgebras}

Let $V = V_0\oplus V_1$ be a superspace of superdimension $(n,n)$, i.e.\
$V$ be a $Z_2$-graded vector space, such that $\dim(V_0)=\dim(V_1)=n$. The set
of linear operators ${\rm End}\, (V)$ acting in the $V$ be an associative $Z_2$-graded
algebra (or superalgebra), ${\rm End}\,(V)=({\rm End}\,(V))_0 \oplus ({\rm End}\,(V))_1$, if the grading
def\/ined by formula:
\[
({\rm End}\,(V))_k=\{g \in {\rm End}\,(V): gV_i \subset V_{i+k}\}.
\]

Let $\{e_1, \ldots, e_n, e_{n+1}, \ldots, e_{2n}\}$
be a such base in $V$ that $\{e_1, \ldots, e_n \}$ be a
basis in $V_0$, $\{e_{n+1}, \ldots$, $e_{2n}\}$ be a basis in $V_1$.
Then we can identify ${\rm End}\,(V)$ with the superalgebra of $(2n \times 2n)$-matrices
$\mathfrak{gl}(n,n)$ (see also \cite{K,F-S}).
Let us def\/ine on the homogeneous components of $\mathfrak{gl}(n,n)$ the commutator
(or supercommutator) by the formula:
\[
[A,B]=AB - (-1)^{\deg(A)\deg(B)}BA,
\]
where $\deg(A)=i$ for $A \in \mathfrak{gl}(n,n)_i$, $i \in Z_2$.
Then $\mathfrak{gl}(n,n)$ turns into the Lie superalgebra.
Further, we will numerate the vectors of the base of $V$ by
integer numbers $\pm 1, \ldots, \pm n$, i.e.\ $\{e_1, \ldots, e_n, e_{-1}, \ldots, e_{-n}\}$ is a basis in $V$,
$\{e_1, \ldots, e_n \}$ is a basis in $V_0$, $\{e_{-1}, \ldots, e_{-n}\}$
 is a basis in $V_1$. Then matrices from $\mathfrak{gl}(n,n)$
 are indexed by numbers: $\pm 1, \ldots, \pm n$, also. Note that
\begin{gather*}
\mathfrak{gl}(n,n)_0= \left\{ \begin{pmatrix} A & 0\\ 0 & B\\ \end{pmatrix}: A, B
\in \mathfrak{gl}(n)\right\},\\
\mathfrak{gl}(n,n)_1= \left\{ \begin{pmatrix} 0 & C\\ D & 0\\ \end{pmatrix}: C, D
\in \mathfrak{gl}(n)\right\}.
\end{gather*}
Let $A=\begin{pmatrix} A_{11} & A_{12}\\ A_{21} & A_{22}\\ \end{pmatrix},
A_{ij} \in \mathfrak{gl}(n)$. By def\/inition,  put:
\[
{\rm str}\,(A)={\rm tr}\,(A_{11})-{\rm tr}\,(A_{22}).
\]
The ${\rm str}\,(A)$ is called the supertrace of $A$.
Def\/ine superalgebra $\mathfrak{sl}(n,n)$ by the formula:
\[
\mathfrak{sl}(n,n)=\{A \in \mathfrak{gl}(n,n): {\rm str}\,(A)=0\}.
\]
Let us also denote $\mathfrak{sl}(n,n)$ by $\tilde{A}(n-1,n-1)$. The Lie
superalgebra $\tilde{A}(n-1,n-1)$ has 1-dimensional center $Z$. Then
$A(n-1, n-1):=\tilde{A}(n-1,n-1)/Z$ is a simple Lie superalgebra.

Let $\pi : \tilde{A}(n-1,n-1) \rightarrow A(n-1,n-1)$ be a natural projection.

Consider the isomorphism  $\sigma': \tilde{A}(n-1, n-1) \rightarrow \tilde{A}(n-1, n-1)$,
which is def\/ined on matrix units $E_{i,j}$ by the formula $\sigma (E_{i,j}) =E_{-i,-j}$.
As  $\sigma'(Z)=Z$, then $\sigma'$ induces the involutive
automorphism  $\sigma : A(n-1, n-1) \rightarrow A(n-1, n-1)$.
Let $\mathfrak{g} = A(n-1, n-1)$. As  $\sigma^2=1$,
then eigenvalues of $\sigma$ equal  $\pm 1$. Let $\epsilon = -1$, $j \in Z_2= \{0,1\}$.
Let us set  $
\mathfrak{g}^j = {\rm Ker}\,(\sigma - \epsilon^jE)$, $\mathfrak{g}=\mathfrak{g}^0 \oplus \mathfrak{g}^1$.
We emphasize that $\mathfrak{g}^0 = \mathfrak{g}^{\sigma}$ is
a set of f\/ixed points of automorphism $\sigma$. It is clear that $\mathfrak{g}^{\sigma}$
is a Lie subsuperalgebra of Lie superalgebra $\mathfrak{g}$. By def\/inition
$Q_{n-1} = \mathfrak{g}^{\sigma}$  is a strange Lie superalgebra. Its inverse
image in  $\tilde{A}(n-1, n-1)$ we denote by $\tilde{Q}_{n-1}$.

We will use the following properties of the Lie superalgebra $Q_{n-1}$.
The root system $\Delta$ of the Lie superalgebra $Q_{n-1}$ coincides
with root system of the Lie algebra $A_{n-1}=\mathfrak{sl}(n)$,
but the non-zero roots of $Q_{n-1}$ are both even and odd.
We will use the following notations: $A = (a_{ij})_{i,j =1}^{n-1}$
is a~Cartan matrix of $\mathfrak{sl}(n)$, $(\alpha_i, \alpha_j)=a_{ij}$
for simple roots $\alpha_i$, $\alpha_j$  ($i,j \in \{1, \ldots,n-1\}$).
Def\/ine the generators of the Lie superalgebra $Q_{n-1}$ $x^{\pm}_i$, $\hat{x}^{\pm}_i$, $h_i$, $k_i$, $i=1, \ldots, n-1$
and elements ${x^{\pm}}{}^i$, ${\hat{x}^{\pm}}{}^i$, $h^i$, $k^i$ of $\mathfrak{g}^1$ by formulas
\begin{gather*}
 h_i = \pi((E_{i,i} -E_{i+1,i+1}) + (E_{i,i} -E_{-i-1,-i-1})), \\
 h^i = \pi((E_{i,i} -E_{i+1,i+1}) - (E_{i,i} -E_{-i-1,-i-1})), \\
x^+_i = \pi(E_{i,i+1}+E_{-i, -i-1}), \qquad {x^+}{}^i = \pi(E_{i,i+1} -E _{-i, -i-1}), \\
x^-_i = \pi(E_{i+1,i}+E_{-i-1, -i}), \qquad {x^-}{}^i = \pi(E_{i+1,i} -E _{-i-1, -i}), \\
 k_i = \pi((E_{i,-i} -E_{i+1,-i-1}) + (E_{-i,i} -E_{-i-1,i+1})),\\
k^i = \pi((E_{i,-i} -E_{i+1,-i-1}) - (E_{-i,i} -E_{-i-1,i+1})), \\
\hat{x}^+_i = \pi(E_{i,-i-1}+E_{-i, i+1}), \qquad {\hat{x}^+}{}^i = \pi(E_{i,-i-1} -E _{-i, i+1}), \\
\hat{x}^-_i = \pi(E_{i+1,-i}+E_{-i-1, i}), \qquad {\hat{x}^-}{}^i = \pi(E_{i+1,-i} -E _{-i-1, i}).
\end{gather*}
    The Lie superalgebra $Q_{n-1}$ can be def\/ined as
    superalgebra generated by generators $h_i$, $k_i$, $x^{\pm}_i$, $\hat{x}^{\pm}_i$, $i \in \{1, \ldots , n-1 \},$
    satisfying the commutation relations of Cartan--Weyl type (see~\cite{F-S}).
    We will use notations $x_{\pm \alpha_i} = x^{\pm}_i$, $\hat{x}_{\pm \alpha_i}
    = \hat{x}^{\pm}_i$, $x^{\pm \alpha_i} = {x^{\pm}}^i$, $\hat{x}^{\pm \alpha_i} = {\hat{x}^{\pm}}{}^i.$
   There exists a nondegenerate supersymmetric invariant bilinear form
   $(\cdot ,\cdot)$ on the Lie superalgebra $A(n-1,n-1)$
   such that $(\mathfrak{g}^0, \mathfrak{g}^0)= (\mathfrak{g}^1, \mathfrak{g}^1)=0$
   and $\mathfrak{g}^0$ and $\mathfrak{g}^1$ nondegenerately paired.
   We use also root genera\-tors~$x_{\alpha}$, $\hat{x}_{\alpha}$ $(\alpha \in \Delta)$
   and elements $x^{\alpha}, \hat{x}^{\alpha} \in \mathfrak{g}^1$ dual (relatively form $(\cdot ,\cdot)$) to them.

Let us extend the automorphism $\sigma$ to automorphism $\tilde{\sigma}:
\mathfrak{g}((u^{-1})) \rightarrow \mathfrak{g}((u^{-1}))$, on
Laurent series with values in $\mathfrak{g}$ by formula:
\begin{gather*}
\tilde{\sigma}(x\cdot u^j)= \sigma(x)(-u)^j.
\end{gather*}
Consider the following Manin triple $(\mathfrak{P}, \mathfrak{P}_1, \mathfrak{P}_2)$:
\[
(\mathfrak{P}=\mathfrak{g}((u^{-1}))^{\tilde{\sigma}},
\qquad
\mathfrak{P}_1= \mathfrak{g}[u]^{\tilde{\sigma}}, \qquad
\mathfrak{P}_2=(u^{-1}\mathfrak{g}[[u^{-1}]])^{\tilde{\sigma}}.
\]

Def\/ine the bilinear form $\langle \cdot, \cdot\rangle$ on $\mathfrak{P}$ by the formula:
\begin{gather} \label{equation20}
\langle f,g\rangle={\rm res}\,(f(u),g(u))du,
\end{gather}
where ${\rm res}\big(\sum\limits_{k=-\infty}^n a_k \cdot u^k\big):= a_{-1}$, $(\cdot, \cdot)$
is an invariant bilinear form on $\mathfrak{g}$.
 It is clear that $\mathfrak{P}_1$, $\mathfrak{P}_2$
 are isotropic subsuperalgebras in relation to the form $\langle\cdot, \cdot\rangle$.
 It can be shown in the usual way that there are the following decompositions:
\begin{gather} \label{equation31}
\mathfrak{g}[u]^{\tilde{\sigma}} =
\bigoplus_{k=0}^{\infty}\big(\mathfrak{g}^0 \cdot u^{2k} \oplus \mathfrak{g}^1 \cdot u^{2k+1}\big) ,
\qquad
\mathfrak{g}((u^{-1}))^{\tilde{\sigma}} =
\bigoplus_{k \in Z}\big(\mathfrak{g}^0 \cdot u^{2k} \oplus \mathfrak{g}^1 \cdot u^{2k+1}\big).
\end{gather}
Describe the bisuperalgebra structures on $\mathfrak{g}[u]^{\tilde{\sigma}}$. Let $\{e_i\}$
 be a basis in $\mathfrak{g}^0$ and  $\{e^i\}$ be a dual basis in $\mathfrak{g}^1$   in relation to
 the form  $(\cdot, \cdot)$. Let $t_0 = \sum e_i \otimes e^i$, $t_1 = \sum e^i \otimes e_i$, $t = t_0 + t_1$.
 Consider also the basis $\{e_{i,k}\}$ in $\mathfrak{P}_1$
and dual basis $\{e^{i,k}\}$ ($\subset \mathfrak{P}_2$) in relation to
the form  $\langle\cdot, \cdot\rangle$, which def\/ine by the formulas:
\begin{gather*}
e_{i,2k}=e_i \cdot u^{2k}, \qquad e_{i,2k+1}=e^i \cdot u^{2k+1}, \qquad k \in Z_+,\\
e^{i,2k}=e^i \cdot u^{-2k-1}, \qquad e^{i,2k+1}=e_i \cdot u^{-2k-2}, \qquad k \in Z_+.
\end{gather*}

Calculate a canonical element $r$, def\/ining the cocommutator in $\mathfrak{P}$
\begin{gather*}
r= \sum e_{i,k} \oplus e^{i,k} =
\sum_{k \in Z} \sum_i \big( e_i \cdot v^{2k} \otimes e^i \cdot u^{-2k-1} +
e^i \cdot v^{2k+1} \otimes e_i \cdot u^{-2k-2}\big) \\
\phantom{r}=\sum_{k=0}^{\infty}\left(\left(\sum e_i \otimes e^i\right) \cdot u^{-1}\left(\dfrac{v}{u}\right)^{2k}\right) +
\sum_{k=0}^{\infty} \left(\left(\sum e^i \otimes e_i\right) \cdot u^{-1} \left(\dfrac{v}{u}\right)^{2k+1}\right) \\
\phantom{r}=t_0 \dfrac{u^{-1}}{(1 - (v/u)^2)} +
t_1 \dfrac{u^{-1}(v/u)}{1 - (v/u)^2} = \dfrac{t_0 \cdot u}{(u^2 - v^2)} +
\dfrac{t_1 \cdot v}{u^2 - v^2} \\
\phantom{r}=  \frac{1}{2} (\dfrac{1}{u-v} + \dfrac{1}{u+v}) t_0 +
\frac 12 \left(\dfrac{1}{u-v} - \dfrac{1}{u+v}\right) t_1 \\
\phantom{r}= \dfrac{1}{2} \dfrac{t_0 + t_1}{u-v} +
\frac{1}{2} \dfrac{t_0 - t_1}{u+v}=
\frac{1}{2} \sum_{k \in Z_+} \dfrac{(\sigma^k \otimes id)\cdot t}{u-\epsilon^k \cdot v}.
\end{gather*}

Denote $r_{\sigma}(u,v):=r$. Then we have the following expression for cocommutator:
\[
\delta : a(u) \rightarrow [a(u) \otimes 1 + 1 \otimes a(v), r_{\sigma}(u,v)].
\]

\begin{proposition} The element $r_{\sigma}(u,v)$ has the following properties:
\begin{gather*}
1) \ \ r_{\sigma}(u,v)= -r_{\sigma}^{21}(v,u);\\
2) \ \ \big[r_{\sigma}^{12}(u,v), r_{\sigma}^{13}(u,w)\big] + \big[r_{\sigma}^{12}(u,v), r_{\sigma}^{23}(v,w)\big]+
\big[r_{\sigma}^{13}(u,w), r_{\sigma}^{23}(v,w)\big] = 0.
\end{gather*}
\end{proposition}

\begin{proof} Let us note, that $t_0^{21}=t_1$, $t_1^{21}=t_0$. Then
\[
r_{\sigma}^{21}(v,u)=
\frac{1}{2} \dfrac{t_1 +t_0}{v-u} + \frac{1}{2} \dfrac{t_1 - t_0}{v+u} = -r_{\sigma}(u,v).
\]
The 1) is proved. Property 2) follows from the fact that the element $r$ def\/ined
above satisf\/ies the classical Yang--Baxter equation. Actually, the function
$r(u,v)= \frac{t_0 + t_1}{u-v}= \frac{t}{u-v}$ satisf\/ies the classical Yang--Baxter equation (CYBE) (see
 \cite{L-S,Dr}):
\begin{gather} \label{equation21}
\big[r^{12}(u,v), r^{13}(u,w)\big] + \big[r^{12}(u,v), r^{23}(v,w)\big]+\big[r^{13}(u,w), r^{23}(v,w)\big] =0.
\end{gather}
Let $s= - {\rm id}$. Let us apply to the left-hand side of (\ref{equation21}) the operator
${\rm id} \otimes s^k \otimes s^l$ $(k,l \in Z_2)$
and substitute $(-1)^k \cdot v$, $(-1)^l \cdot w$ for $v$, $w$, respectively.
Taking then the sum over $k,l \in Z_2$ and using  $(s \otimes s)(t)  = -t$
we will obtain the left-hand side of expression from item 2 of proposition.
\end{proof}

\section{Quantization}

The def\/inition of quantization of Lie bialgebras from \cite{Dr}
can be naturally extended on Lie bisuperalgebras.
Quantization of Lie bisuperalgebras $D$
is such Hopf superalgebra $A_\hbar$ over ring of formal power
series $\cal{C}[[\hbar]]$ that satisfies the following conditions:
\begin{enumerate}\itemsep=0pt
\item[1)] $A_\hbar/\hbar A_\hbar \cong U(D),$
as a Hopf algebra (where $U(D)$  is an universal enveloping algebra of the Lie superalgebra $D$);
\item[2)] the superalgebra $A_\hbar$ isomorphic to $U(D)[[\hbar]],$ as a vector space;
\item[3)] it is fulf\/illed the following correspondence principle: for any $x_0 \in D$
and any $x \in A_\hbar$ equal to $x_0$:   $x_0 \equiv x \mod \hbar$ one has
\[
\hbar^{-1}(\Delta(x)-\Delta^{\rm op}(x)) \mod \hbar \equiv  \varphi(x) \mod \hbar,
\]
 where $\Delta$ is a comultiplication, $\Delta^{\rm op}$ is an opposite comultiplication
 (i.e., if  $\Delta(x)= \sum x'_i\otimes x''_i$, then $\Delta^{\rm op}(x)=\sum(-1)^{p(x'_i)p(x''_i)}x''_i\otimes x'_i$).
\end{enumerate}

Let us describe the quantization of Lie bisuperalgebra
$(\mathfrak{g}[u]^{\tilde{\sigma}}, \delta)$. I recall (see (\ref{equation31})) that
\[
\mathfrak{g}[u]^{\tilde{\sigma}} =
\bigoplus_{k=0}^{\infty}
\big(\mathfrak{g}^0 \cdot u^{2k} \oplus \mathfrak{g}^1 \cdot u^{2k+1}\big)
\] is graded by degrees of $u$  Lie superalgebra,
\begin{gather} \label{equation30}
\delta : a(u) \rightarrow \left[a(u)\otimes 1 + 1 \otimes a(v),
\frac{1}{2} \sum_{k \in Z_+} \dfrac{(\sigma^k \otimes {\rm id})\cdot \mathfrak{t}}{u-\epsilon^k \cdot v}\right],
\end{gather}
where $\mathfrak{t}$ is a Casimir operator and  $\delta$ is a  homogeneous map of degree $-1$.

 Let us apply  the additional conditions upon quantization.
 \begin{enumerate}
\item[1)] Let $A$ be a graded superalgebra over graded ring $C[[\hbar]]$, $\deg(\hbar)=1$.
\item[2)] The grading of $A$ and the grading of  \(\mathfrak{g}[u]^{\tilde{\sigma}}\)
induce the same gradings of \(U(\mathfrak{g}[u]^{\tilde{\sigma}})\), i.e.
\[
A/\hbar A = U(\mathfrak{g}[u]^{\tilde{\sigma}})
\] as graded superalgebra over $\mathbb{C}$.
\end{enumerate}

I recall (see also \cite{Dr}) that Hopf superalgebra $A$ over $C[[\hbar]]$
such that $A/\hbar A \cong B$, where $B$ is a Hopf superalgebra
over $\mathbb{C}$, is called a formal deformation of $B$.
Let \(p: A \rightarrow A/\hbar A \cong B \)  be a canonical projection.
If  $p(a)=x,$ then an element $a$ is called a deformation of element~$x$.
There exist theorems that prove existence and uniqueness of
quantization  (or formal deformation) in many special cases as well as in our one.
But, we will not use these theorems. Let $m_i \in \{h_i, k_i, x^{\pm}_i, \hat{x}^{\pm}_i\}$
be generators of the Lie superalgebra $Q_n$. We will denote
the deformations of the generators $m_i\cdot u^k$ of the
Lie superalgebra $A(n,n)[u]^{\tilde{\sigma}}$ by $m_{i,k}$.
As generators of associative superalgebra $m_i$, $m_i \cdot u^k$ generate
the superalgebra $U(A(n,n)[u]^{\tilde{\sigma}})$, its deformations $m_{i,0}$, $m_{i,1}$
generate the Hopf superalgebra $A$. We are going to describe
the system of def\/ining relations between these generators.
This system of def\/ining relations is received  from conditions of compatibility
of superalgebra and cosuperalgebra  structures of $A$ (or
from condition that comultiplication is a homomorphism of superalgebras).
First, we describe the comultiplication on generators $m_{i,1}$.
It follows from  condition of homogeneity of quantization 2) comultiplication
is def\/ined only by of values of $\Delta$ on generators $m_{i,1}$.
Describe values of $\Delta$ on generators $h_{i,1}$. Note
that condition of homogeneity of quantization implies the fact that
\(U(\mathfrak{g}^0)\) embeds in $A$ as a Hopf superalgebra.
It means that  we can identify generators $m_{i,0}$ with generators
of Lie superalgebra $\mathfrak{g}^0=Q_{n-1}$. Calculate the value of
cocycle $\delta$ on the generators $h_i \cdot u$, $i=1, \dots , n-1$.

 \begin{proposition}
\label{proposition1}
Let $\mathfrak{A}$ be a Lie superalgebra with
invariant scalar product $(\cdot, \cdot)$; $\{e_i\}$, $\{e^i\}$ be a~dual
relatively this scalar product bases. Then for every element $g \in \mathfrak{A}$ we have equality:
\begin{gather} \label{equation311}
\left[g \otimes 1, \sum e_i \otimes e^i\right] = -\left[1 \otimes g, \sum e_i \otimes e^i\right] .
\end{gather}
\end{proposition}

\begin{proof} Note that following the def\/inition of bilinear invariant form, we have the
equality: \( ([g, a], b)= -(-1)^{\deg(g)\deg(a)}([a,g],b)
= -(-1)^{\deg(g)\deg(a)}(a,[g,b])\) for $\forall \,a, b \in \mathfrak{A}$.  Therefore
\begin{gather}\label{equation32}
([g, e_i], e^i)= -(-1)^{\deg(g)\deg(a)}([e_i,g],e^i)= -(-1)^{\deg(g)\deg(a)}(e_i,[g,e^i]).
\end{gather}
The scalar product given on vector space $V$ def\/ines isomorphism between $V$ and $V^*$,
and, therefore, between $V\otimes V$ and $V\otimes V^*$. Summing over $i$
equality (\ref{equation32}) we get
\[ \sum_i ([g, e_i], e^i)= \sum_i -(-1)^{\deg(g)\deg(a)}(e_i,[g,e^i]). \]
Note the equality of functionals follows from the equality of values of
functionals on elements of base and thus we have an equality:
\[
\sum_i [g, e_i]\otimes e^i= \sum_i -(-1)^{\deg(g)\deg(a)}e_i \otimes [g,e^i]
 \]
or  $\big[g \otimes 1, \sum e_i \otimes e^i\big] = -\big[1 \otimes g, \sum e_i \otimes e^i\big]$ or equality (\ref{equation311}).
\end{proof}

\begin{proposition}
Let \(\mathfrak{A}= \mathfrak{A}^0 \oplus \mathfrak{A}^1\)
be a Lie superalgebra with such nondegenerate invariant scalar
product that $\mathfrak{A}^0$, $\mathfrak{A}^1$ are isotropic subspaces, $\mathfrak{A}^0$, $\mathfrak{A}^1$
are nondegenerately paired, $\mathfrak{A}^0$  is a~subsuperalgebra, $\mathfrak{A}^1$ is a module over
\(\mathfrak{A}^0\). (For example, the scalar product \eqref{equation20} satisfies
these conditions.)
Let also $\{e_i\}$, $\{e^i\}$ be the dual bases
in  $\mathfrak{A}^0$, $\mathfrak{A}^1$, respectively,
and $\mathfrak{t}_0= \sum_i e_i \otimes e^i$, $\mathfrak{t}_1= \sum_i e^i \otimes e_i$.
 Then for all  $a \in \mathfrak{A}^0$, $b \in \mathfrak{A}^1$  we have the following equalities:
\begin{gather*}
[a\otimes 1, \mathfrak{t}_0] = - [1 \otimes a, \mathfrak{t}_0],\qquad
[a\otimes 1, \mathfrak{t}_1] = - [1 \otimes a, \mathfrak{t}_1],\\
[b\otimes 1, \mathfrak{t}_0] = - [1 \otimes b, \mathfrak{t}_1],\qquad
[b\otimes 1, \mathfrak{t}_1] = - [1 \otimes b, \mathfrak{t}_0].
\end{gather*}
\end{proposition}

 Now we can calculate the value $\delta$ (see (\ref{equation30})) on $h^i\cdot u$
\begin{gather*}
 \delta(h^i \cdot u) = \left[h^i \cdot v \otimes 1 + 1 \otimes h^i \cdot u, \frac{1}{2} \dfrac{\mathfrak{t}_0 + \mathfrak{t}_1}{u - v} + \frac{1}{2} \dfrac{\mathfrak{t}_0 - \mathfrak{t}_1}{u + v}\right] \\
\phantom{\delta(h^i \cdot u)}{} =\left[h^i \cdot v \otimes 1 - h^i \cdot u \otimes 1, \frac{1}{2} \dfrac{\mathfrak{t}_0 + \mathfrak{t}_1}{u - v}\right] + \left[ h^i \cdot v \otimes 1 + h^i\cdot u \otimes 1, \frac{1}{2} \dfrac{\mathfrak{t}_0 - \mathfrak{t}_1}{u + v} \right] \\
\phantom{\delta(h^i \cdot u)}{}=   \left[h^i \otimes 1, \frac{1}{2} (\mathfrak{t}_0 + \mathfrak{t}_1)\right] + \left[h^i \otimes 1, \frac{1}{2}(\mathfrak{t}_0 - \mathfrak{t}_1)\right]=[h^i \otimes 1, \mathfrak{t}_0]= -[ 1 \otimes h^i, \mathfrak{t}_1].
\end{gather*}
Similarly, it is possible to calculate the values of cocycle on other generators
\begin{gather*}
 \delta(k^i \cdot u) = -[k^i \otimes 1, \mathfrak{t}_0]= [1 \otimes k^i, \mathfrak{t}_1],\\
 \delta({x^{\pm}}{}^i \cdot u) = -[{x^{\pm}}{}^i \otimes 1, \mathfrak{t}_0]= [1 \otimes {x^{\pm}}{}^i, \mathfrak{t}_1],\\
 \delta({\hat{x}^{\pm}}{}^i \cdot u) = -[{\hat{x}^{\pm}}{}^i \otimes 1, \mathfrak{t}_0]
= [1 \otimes {\hat{x}^{\pm}}{}^i, \mathfrak{t}_1].
\end{gather*}

  It follows from homogeneity condition that
\begin{gather*}
\Delta(h_{i,1})= \Delta_0(h_{i,1}) + \hbar F(x_{\alpha}\otimes x_{-\alpha}, \hat{x}_{\alpha} \otimes \hat{x}_{-\alpha} + h_i \otimes h_j + k_i \otimes k_j).
\end{gather*}
It follows from correspondence principle (item 3) of def\/inition of quantization) that
\[
\hbar^{-1}(\Delta(h_{i,1})- \Delta^{\rm op}(h_{i,1}) = F - \tau F = [1 \otimes h^i, \mathfrak{t}_0].
\]

Let {\samepage
\[
\bar{\mathfrak{t}}_0 = \sum_{\alpha \in \Delta_+} x_{\alpha} \otimes x^{-\alpha} - \hat{x}_{\alpha} \otimes \hat{x}^{-\alpha} + \frac{1}{2} \sum_{i=1}^{n-1} k_i \otimes k^i,
\]
\(\Delta_+ \) is a set of positive roots of Lie algebra  $A_{n-1}= \mathfrak{sl}(n)$.}

Def\/ine \(\Delta(h_{i,1})\),
by formula
\[
\Delta(h_{i,1})= \Delta_0(h_{i,1}) + \hbar [1 \otimes h^i, \bar{\mathfrak{t}}_0].
\]
Let us check that correspondence principle is fulf\/illed in these cases too:
\begin{gather*}
\hbar^{-1}(\Delta(h_{i,1})- \Delta^{\rm op}(h_{i,1})) =
\left[1 \otimes h^i, \sum_{\alpha \in \Delta_+} x_{\alpha} \otimes x^{-\alpha} - \hat{x}_{\alpha} \otimes \hat{x}^{-\alpha} + \frac{1}{2} \sum_{i=1}^{n-1} k_i \otimes k^i\right]\\
 \qquad{}- \left[h^i \otimes 1,  \sum_{\alpha \in \Delta_+} x^{-\alpha} \otimes x_{\alpha} - \hat{x}^{-\alpha} \otimes \hat{x}_{\alpha} + \frac{1}{2} \sum_{i=1}^{n-1} k^i \otimes k_i\right]\\
  \qquad{}= \left[1 \otimes h^i, \sum_{\alpha \in \Delta_+} x_{\alpha} \otimes x^{-\alpha} - \hat{x}_{\alpha} \otimes \hat{x}^{-\alpha}\right] - \left[h^i \otimes 1, \sum_{\alpha \in \Delta_+} x^{-\alpha} \otimes x_{\alpha} - \hat{x}^{-\alpha} \otimes \hat{x}_{\alpha}\right].
\end{gather*}

Show that
\[
\left[h^i \otimes 1, \sum x^{-\alpha}\otimes x_{\alpha} + \hat{x}^{-\alpha}\otimes \hat{x}_{\alpha}\right] =
\left[1 \otimes h^i, \sum x_{-\alpha}\otimes x^{\alpha} + \hat{x}_{-\alpha}\otimes \hat{x}^{\alpha}\right].
\]
Actually,
\[
\left[h^i \otimes 1, \sum x^{-\alpha}\otimes x_{\alpha} + \hat{x}^{-\alpha}\otimes \hat{x}_{\alpha}\right] = \sum [h^i, x^{-\alpha}]\otimes x_{\alpha} + [h^i, \hat{x}^{-\alpha}] \otimes \hat{x}_{\alpha}.
\]
On the other hand
\[
\left[1 \otimes h^i, \sum x_{-\alpha}\otimes x^{\alpha} + \hat{x}_{-\alpha}\otimes \hat{x}^{\alpha}\right] = \sum x_{-\alpha}\otimes [h^i, x^{\alpha}] + \hat{x}_{-\alpha} \otimes [h^i, \hat{x}^{\alpha}].
\]

It follows in the standard way, that
\begin{gather*}
[h^i, x^{\alpha_i - \alpha_j}] = -(\delta_{ik} + \delta_{jk} - \delta_{i,k+1} - \delta_{j,k+1})x_{\alpha_i - \alpha_j},\\
[k^i, \hat{x}^{\alpha_i - \alpha_j}] = (\delta_{ik} - \delta_{jk} - \delta_{i,k+1} + \delta_{j,k+1})\hat{x}_{\alpha_i - \alpha_j}.
\end{gather*}
Therefore
\begin{gather*}
[h^i, x^{-\alpha}]\otimes x_{\alpha}= - x_{-\alpha} \otimes [h^i, x^{\alpha}],\qquad
[h^i, \hat{x}^{-\alpha}]\otimes \hat{x}_{\alpha}= \hat{x}_{-\alpha} \otimes [h^i,\hat{x}^{\alpha}].
\end{gather*}
Hence,
\begin{gather*}
\left[h^i \otimes 1, \sum x^{-\alpha}\otimes x_{\alpha} + \hat{x}^{-\alpha}\otimes \hat{x}_{\alpha}\right] = \left[1 \otimes h^i, \sum x_{-\alpha}\otimes x^{\alpha} + \hat{x}_{-\alpha}\otimes \hat{x}^{\alpha}\right]
\end{gather*}
and equality
\begin{gather*}
\hbar^{-1}(\Delta(h_{i,1})- \Delta^{\rm op}(h_{i,1})) = \left[1 \otimes h^i, \sum_{\alpha \in \Delta_+} x_{\alpha} \otimes x^{-\alpha} -\hat{x}_{\alpha} \otimes \hat{x}^{-\alpha}\right]\\
\qquad{} -\left[h^i \otimes 1, \sum_{\alpha \in \Delta_+} x^{-\alpha} \otimes x_{\alpha} - \hat{x}^{-\alpha} \otimes \hat{x}_{\alpha}\right]
= [1\otimes h^i, \mathfrak{t}_0]= \delta(h^i\cdot u)
\end{gather*}
is proved.

Similarly we can def\/ine comultiplication on other generators. We get
\begin{gather*}
 \Delta(x^+_{i,1})=\Delta_0(x^+_{i,1}) + \hbar [1\otimes x^+_i, \bar{\mathfrak{t}}_0], \qquad
\Delta(x^-_{i,1})=\Delta_0(x^-_{i,1}) + \hbar [x^-_i \otimes 1, \bar{\mathfrak{t}}_0].
\end{gather*}
The following relations hold and comultiplication preserves them
\begin{gather*}
 [h_{i,1}, x^{\pm}_{j,0}]= \pm (\alpha_, \alpha_j)x^{\pm}_{j,1}, \qquad
[k_{i,1}, x^{\pm}_{j,0}]= \pm (\alpha_, \alpha_j)\hat{x}^{\pm}_{j,1}, \\
[h_{i,1}, \hat{x}^{\pm}_{j,0}]= \pm \widetilde{(\alpha_, \alpha_j)}\hat{x}^{\pm}_{j,1}, \qquad
[k_{i,1}, \hat{x}^{\pm}_{j,0}]= \pm (\alpha_, \alpha_j)x^{\pm}_{j,1}.
\end{gather*}

   Let us prove  these relations. It is sufficient to check one of them. Let us prove the f\/irst relation
\[
[\Delta(h_{i,1}), \Delta(x^{\pm}_{j,0})]= \pm (\alpha_i, \alpha_j)\Delta(x^{\pm}_{j,1}).
\]
Actually,
\begin{gather*}
[\Delta(h_{i,1}), \Delta(x^{\pm}_{j,0})]= [\Delta_0(h_{i,1}) + \hbar [1\otimes h_i, \bar{\mathfrak{t}}_0] , x^{\pm}_{j,0} \otimes 1 + 1\otimes x^{\pm}_{j,0} ]\\
\qquad{}= [h_{i,1}, x^{\pm}_{j,0}]\otimes 1 + 1 \otimes [h_{i,1}, x^{\pm}_{j,0}] + \hbar [[1 \otimes h_i, \bar{\mathfrak{t}}_0], x^{\pm}_{j,0} \otimes 1] + \hbar[[1 \otimes h_i, \bar{\mathfrak{t}}_0], 1 \otimes x^{\pm}_{j,0}]\\
\qquad{} = [\Delta(h_{i,0}), \Delta(x^{\pm}_{j,1})] =  \pm (\alpha_i, \alpha_j)\Delta(x^{\pm}_{j,1}).
\end{gather*}

 Now we can describe the Hopf superalgebra $A=A_h$, which is a deformation (quantization) of the Lie bisuperalgebra $(\mathfrak{g}^{\tilde{\sigma}}, \delta)$. Let us introduce new notations.
Let $\{a,b\}$ be an anticommutator of elements $a$, $b$.  Let  $m \in \{0,1\}$ and
\begin{gather*}
\bar {k}_{i,m} = {1 \over n}
\left( - \sum\limits_{r = 0}^{i - 1}{rk_{r,m} } + \sum\limits_{r = i}^{n - 1} {(n - r)k_{r,m} } \right), \\
\bar {h}_{i,m} = {1 \over n}\left( - \sum\limits_{r = 0}^{i - 1} {rh_{r,m} } +
\sum\limits_{r = i}^{n - 1} {(n - r)h_{r,m} } \right).
\end{gather*}
Let also, $\{\alpha_1, \cdots, \alpha_{n-1}\}$ be a set of simple roots of $\mathfrak{sl}(n)$,
$\widetilde{(\alpha _i ,\alpha _j)} : = (\delta _{i,j + 1} - \delta _{i + 1,j})$.

\begin{theorem}\label{theorem1}
Hopf superalgebra $A=A_h$ over $\mathcal{C}[[\hbar]]$ is generated by generators $h_{i,0}$, $x^{\pm}_{i,0}$, $k_{i,0}$, $\hat{x}^{\pm}_{i,0}$,
$h_{i,1}$, $x^{\pm}_{i,1}$, $k_{i,1}$, $\hat{x}^{\pm}_{i,1}$, $1\leq i\leq n-1$ ($h$, $x$ are even, $k$, $\hat{x}$ are odd generators). These generators satisfy the following defining relations:
\begin{gather*}
[h_{i,0},h_{j,0}]=[h_{i,0},h_{j,1}]=[h_{i,1},h_{j,1}]=0,\qquad
[h_{i,0},k_{j,0}]=[k_{i,1},k_{j,0}]=0,  \\
[h_{i,0},x_{j,0}^{\pm}]=\pm (\alpha_i,\alpha_j) x_{j,0}^{\pm}, \qquad
[k_{i,0},x_{j,0}^{\pm}]=\pm \widetilde{(\alpha_i,\alpha_j)} \hat{x}_{j,0}^{\pm}, \\
[k_{i,0}, k_{j,0}]=2(\delta_{i,j}-\delta_{i,j+1})\bar{h}_{i,0} + 2(\delta_{i,j}-\delta_{i,j-1})\bar{h}_{i+1,0}, \\
k_{i,1}=\frac{1}{2}[h_{i+1,1}-h_{i-1,1}, k_{i,0}], \qquad
[x_{i,0}^+, x_{j,0}^-]=\delta_{ij} h_{i,0},  \\
[x_{i,1}^+,x_{j,0}^-]= [x_{i,0}^+,x_{j,1}^-] = \delta_{ij}\tilde{h}_{i,1} =  \delta_{ij}\left(h_{i,1}+\frac{\hbar}{2}h_{i,0}^2\right),\\ 
[\hat{x}_{i,0}^+, x_{j,0}^-]=[x_{i,0}^+, \hat{x}_{j,0}^-]=\delta_{ij} k_{i,0}, \qquad
[x_{i,1}^+, x_{j,0}^-]=\delta_{ij} \left(h_{i,1} +\frac{\hbar}{2}h_{i,0}^2\right), \\
[\hat{x}_{i,1}^+, x_{j,0}^-]=[x_{i,0}^+, \hat{x}_{j,1}^-]=\delta_{ij} k_{i,1}, \qquad
[x_{i,1}^+, \hat{x}_{j,0}^-]=-[\hat{x}_{i,0}^+, x_{j,1}^-]=\delta_{ij} (\bar{k}_{i,1}+ \bar{k}_{i+1,1}),  \\
[\hat{x}_{i,1}^+, \hat{x}_{j,0}^-]=\delta_{ij} \left(h_{i,1} +\frac{\hbar}{2}h_{i,0}^2\right), \qquad
[h_{i,1},x_{j,0}^{\pm}] = \pm (\alpha_i,\alpha_j)(x_{j,1}^{\pm}), \\
[h_{i,1}, \hat{x}_{j,0}^{\pm}]=\pm \widetilde{(\alpha _i ,\alpha _j)}(\hat{x}_{j,1}^{\pm}), \qquad
k_{i,1}= \frac{1}{2}[h_{i+1,1} - h_{i-1,1}, k_{i,0}], \\
[k_{i,1},x_{j,0}^{\pm}]=\pm (\alpha_i,\alpha_j)\hat{x}_{j,1}^{\pm}, \qquad 
[k_{i,1}, \hat{x}_{j,0}^{\pm}]=\pm (\alpha_i,\alpha_j)x_{j,1}^{\pm}, \\ 
[h_{i,1}, k_{j,0}]= 2(\delta_{i,j}-\delta_{i,j+1})\bar{k}_{i,1} + 2(\delta_{i,j}-\delta_{i,j-1})\bar{k}_{i+1,1}, \\
[x_{i,1}^{\pm},x_{j,0}^{\pm}]-[x_{i,0}^{\pm},x_{j,1}^{\pm}] = \pm
\frac{\hbar}{2}( (\alpha_i, \alpha_j)\{x_{i,0}^{\pm}, x_{j,0}^{\pm}\} + \widetilde{(\alpha_i, \alpha_j)}\{\hat{x}_{i,0}^{\pm}, \hat{x}_{j,0}^{\pm}\}), \\ 
[\hat{x}_{i,1}^{\pm},x_{j,0}^{\pm}]-[\hat{x}_{i,0}^{\pm},x_{j,1}^{\pm}] = \pm \frac{\hbar}{2}(-\widetilde{(\alpha_i, \alpha_j)}\{\hat{x}_{i,0}^{\pm}, x_{j,0}^{\pm}\} + \widetilde{(\alpha_i, \alpha_j)}\{x_{i,0}^{\pm}, \hat{x}_{j,0}^{\pm}\}),\\ 
[\hat{x}_{i,1}^{\pm}, \hat{x}_{j,0}^{\pm}]-[\hat{x}_{i,0}^{\pm}, \hat{x}_{j,1}^{\pm}] = \pm
\frac{\hbar}{2}(\widetilde{(\alpha_i, \alpha_j)}\{x_{i,0}^{\pm}, x_{j,0}^{\pm}\} + (\alpha_i, \alpha_j)\{\hat{x}_{i,0}^{\pm}, \hat{x}_{j,0}^{\pm}\}), \\ 
({\rm ad}\, x_{i,0}^{\pm})^{2}(x_{j,0}^{\pm})= [x_{i,0}^{\pm},[x_{i,0}^{\pm},x_{j,0}^{\pm}]]= 0,\qquad i\neq j, \\ 
({\rm ad}\,\hat{x}_{i,0}^{\pm})^{2}(x_{j,0}^{\pm})= [\hat{x}_{i,0}^{\pm},[\hat{x}_{i,0}^{\pm}, x_{j,0}^{\pm}]]= 0 = [\hat{x}_{i,0}^{\pm},[x_{i,0}^{\pm}, x_{j,0}^{\pm}]],\qquad i\neq j, \label{equation2351}  \\
({\rm ad}\, x_{i,0}^{\pm})^{2}(\hat{x}_{j,0}^{\pm})= [x_{i,0}^{\pm},[x_{i,0}^{\pm}, \hat{x}_{j,0}^{\pm}]]= 0,\qquad i\neq j, \\ 
[\hat{x}_{i,0}^{\pm}, \hat{x}_{j,0}^{\pm}]= [x_{i,0}^{\pm}, x_{j,0}^{\pm}], \qquad
[\hat{x}_{i,0}^{\pm}, x_{j,0}^{\pm}]= [\hat{x}_{i,0}^{\pm}, x_{j,0}^{\pm}], \\
\sum_{\sigma \in S_3}[x^{\pm}_{i, \sigma(s_1)}, [x^{\pm}_{i, \sigma(s_2)}, x^{\pm}_{i, \sigma(s_3)}]] =0, \\
\sum_{\sigma \in S_3}[\hat{x}^{\pm}_{i, \sigma(s_1)}, [x^{\pm}_{i, \sigma(s_2)}, x^{\pm}_{i, \sigma(s_3)}]] =0, \qquad s_1, s_2, s_3 \in \{0, 1\}.
\end{gather*}
Here $S_n$ is a permutation group of $n$ elements.

The comultiplication $\Delta$ is defined by the formulas:
\begin{gather*}
\Delta (h_{i,0})=h_{i,0}\otimes 1 + 1\otimes h_{i,0},\qquad
\Delta (x_{i,0}^{\pm})=x_{i,0}^{\pm} \otimes 1 + 1\otimes x_{i,0}^{\pm}, \\
\Delta (k_{i,0})=k_{i,0}\otimes 1 - 1\otimes k_{i,0},\qquad
\Delta (\hat{x}_{i,0}^{\pm})=\hat{x}_{i,0}^{\pm} \otimes 1 + 1\otimes \hat{x}_{i,0}^{\pm}, \\
 \Delta (h_{i,1})=h_{i,1}\otimes 1+ 1\otimes h_{i,1}+  [1 \otimes h_{i,0}, \bar{\mathfrak{t}}_0]=h_{i,1}\otimes 1+ 1\otimes h_{i,1}  \\
 \phantom{\Delta (h_{i,1})=}{}+ \hbar\left(k_{i,0} \otimes (\bar{k}_{i,0} +\bar{k}_{i+1,0}) - \sum_{\alpha \in \Delta_+} \big((\alpha_i,\alpha)x_{\alpha, 0} \otimes x_{-\alpha, 0} + \widetilde{(\alpha_i, \alpha)}\hat{x}_{\alpha,0}\otimes \hat{x}_{-\alpha,0}\big)\right), \\
 \Delta(x_{i,1}^+)= x_{i,1}^+\otimes 1+ 1\otimes x_{i,1}^+  +
[1 \otimes x_{i,0}^+,\bar{\mathfrak{t}}_0]=
  x^+_{i,1}\otimes 1+  1\otimes x^+_{i,1}\\
  \phantom{\Delta(x_{i,1}^+)=}{}+ \hbar\Bigg(\hat{x}^+_{i,0} \otimes (\bar{k}_{i,0} + \bar{k}_{i+1,0})  + x^+_{i,0} \otimes h_{i,0}\\
  \phantom{\Delta(x_{i,1}^+)=}{} -
\sum_{\alpha \in \Delta_+} \big([x^+_{i,0}, x_{\alpha,0}] \otimes x_{-\alpha,0} + [x^+_{i,0}, \hat{x}_{-\alpha,0}] \otimes \hat{x}_{-\alpha,0}\big)\Bigg), \\
 \Delta(x_{i,1}^-)= x_{i,1}^-\otimes 1+ 1\otimes x_{i,1}^-  + [x_{i,0}^- \otimes 1,\bar{\mathfrak{t}}_0] \\
\phantom{\Delta(x_{i,1}^-)}{}=  x^-_{i,1}\otimes 1+  1\otimes x^-_{i,1}+ \hbar\Bigg(\big(\bar{k}_{i,0} + \bar{k}_{i+1,0}) \otimes \hat{x}^-_{i,0} + h_{i,0} \otimes x^-_{i,0} \\
\phantom{\Delta(x_{i,1}^-)=}{} -
\sum_{\alpha \in \Delta_+} \big(x_{\alpha,0} \otimes  [x^-_{i,0}, x_{-\alpha,0}] + \hat{x}_{\alpha,0} \otimes [x^-_{i,0}, \hat{x}_{-\alpha,0}]\big)\Bigg), \\
 \Delta(\hat{x}_{i,1}^+)= \hat{x}_{i,1}^+\otimes 1 - 1\otimes \hat{x}_{i,1}^+  +
 \hbar\Bigg(x^+_{i,0} \otimes (\bar{k}_{i,0} + \bar{k}_{i+1,0})  + \hat{x}^+_{i,0} \otimes h_{i,0} \\
 \phantom{\Delta(\hat{x}_{i,1}^+)=} {}-
\sum_{\alpha \in \Delta_+} \big([\hat{x}^+_{i,0}, x_{\alpha,0}] \otimes x_{-\alpha,0} - [x^+_{i,0}, x_{-\alpha,0}] \otimes \hat{x}_{-\alpha,0}\big)\Bigg), \\
 \Delta(\hat{x}_{i,1}^-)= \hat{x}_{i,1}^-\otimes 1 - 1\otimes \hat{x}_{i,1}^-  + [\hat{x}_{i,0}^- \otimes 1,\bar{\mathfrak{t}}_0] \\
 \phantom{\Delta(\hat{x}_{i,1}^-)}{}= \hat{x}^-_{i,1}\otimes 1 +  1 \otimes \hat{x}^-_{i,1}+ \hbar\Bigg((\bar{k}_{i,0} + \bar{k}_{i+1,0}) \otimes x^-_{i,0} + h_{i,0} \otimes \hat{x}^-_{i,0}   \\
\phantom{\Delta(\hat{x}_{i,1}^-)=}{}- \sum_{\alpha \in \Delta_+} \big(x_{\alpha,0} \otimes  [\hat{x}^-_{i,0}, x_{-\alpha,0}] + \hat{x}_{\alpha,0} \otimes [x^-_{i,0}, x_{-\alpha,0}]\big)\Bigg), \\
 \Delta (k_{i,1})=k_{i,1}\otimes 1 - 1\otimes k_{i,1}+  [k_{i,0}\otimes 1,\bar{\mathfrak{t}}_0]  \\
\phantom{\Delta (k_{i,1})}{} = k_{i,1}\otimes 1+ 1\otimes k_{i,1}+ \hbar \Bigg((\bar{h}_{i,0} +\bar{h}_{i+1,0}) \otimes (\bar{k}_{i,0} +\bar{k}_{i+1,0}) \\
\phantom{\Delta (k_{i,1})=}{}-  \sum_{\alpha \in \Delta_+} \big((\alpha_i,\alpha)\hat{x}_{\alpha,0} \otimes x_{-\alpha,0} + \widetilde{(\alpha_i, \alpha)}x_{\alpha,0}\otimes \hat{x}_{-\alpha,0}\big)\Bigg).
\end{gather*}
\end{theorem}

 Let us note that the Hopf superalgebras $A_{\hbar_1}$ and $A_{\hbar_2}$ for f\/ixed $\hbar_1, \hbar_2 \neq 0$ (as superalgebras over~$\mathbb{C}$) are isomorphic. Setting $\hbar=1$ in these formulas we receive the system of def\/ining relations of Yangian  $Y(Q_{n-1})$.

\section{Current system of generators}

Let $G=Q_n$. Let us introduce a new system of  generators and def\/ining relations. This system in the quasiclassical limit transforms to the current system of generators for twisted current Lie superalgebra $\wp _1 = G[u]^{\tilde{\sigma}}$ of polynomial currents. We introduce the generators $\tilde {h}_{i,m}$,
$k_{i,m}$, $x^\pm _{i,m}$, $\tilde {x}^\pm _{i,m}$,
$i \in I = \{1,2,\ldots,n - 1\}$, $m \in Z_ + $, by the following formulas:
\begin{gather}
x^\pm _{i,m + 1} = \pm \frac{1}{2} [h_{i,1} ,x^\pm _{i,m} ], \label{eq37} \\
\hat {x}^\pm _{i,2m + 1} = \frac{1}{2} [h_{i + 1,1}- h_{i - 1,1}, \hat{x}^\pm _{i,2m}], \label{eq38} \\
\hat{x}^\pm _{i,2m + 2} = - \frac{1}{2} [h_{i + 1,1} ,- h_{i-1,1} ,\hat{x}^\pm _{i,2m + 1} ], \label{eq39} \\
k_{i,m + 1} = \frac{1}{2} [h_{i + 1,1} , - h_{i - 1,1} ,k_{i,m} ], \label{eq40} \\
\tilde {h}_{i,m} = [x^ + _{i,m} ,x^ - _{i,0} ],\label{eq41}\\
\bar {k}_{i,m} = {1 \over n}\left( - \sum\limits_{r = 0}^{i - 1}{rk_{r,m} } + \sum\limits_{r = i}^{n - 1} {(n - r)k_{r,m} } \right), \nonumber\\ \bar {h}_{i,m}
= {1 \over n}\left( - \sum\limits_{r = 0}^{i - 1} {r\tilde {h}_{r,m} } +
\sum\limits_{r = i}^{n - 1} {(n - r)\tilde {h}_{r,m} } \right). \label{eq42}
\end{gather}

This section results in to the following theorem describing the $Y(Q_{n - 1})$ in a convenient form.

\begin{theorem}
The Yangian $Y(Q_{n - 1})$ isomorphic to the unital associative superalgebra over $\mathbb C$, generated by generators $\tilde {h}_{i,m}$, $k_{i,m}$, $x^\pm _{i,m}$, $\hat{x}^\pm _{i,m}$, $i \in I = \{1,2,\ldots,n - 1\}$, $m \in Z_ + $  (isomorphism is given by the formulas \eqref{eq37}--\eqref{eq42}), satisfying the following system of defining relations:
\begin{gather*}
[\tilde {h}_{i,m} ,\tilde {h}_{j,n} ] = 0, \qquad
\tilde {h}_{i,m + n} = \delta _{ij} [x^ + _{i,m} ,x^ - _{i,n} ], \label{eq44}\\
[\hat{x}^ + _{i,m} ,x^ - _{j,2k} ] = [x^ + _{i,2k} ,\tilde {x}^ - _{j,m}] =
\delta _{ij} k_{i,m + 2k} \left({{n - 2} \over n}\right)^k, \label{eq45} \\
[\hat{x}^ + _{i,m} ,x^ - _{j,2k + 1} ] = [x^ + _{i,2k + 1} ,\hat{x}^-_{j,m}] = \delta _{ij} (\bar {k}_{i,m + 2k + 1} + \bar {k}_{i + 1,m + 2k + 1} )\left({{n - 2} \over n}\right)^k, \label{eq46}\\
[h_{i,0} ,x^\pm _{j,l} ] = \pm (\alpha _i ,\alpha _j )x^\pm _{j,l} ,\qquad [h_{i,0},
\tilde {x}^\pm _{j,l} ] = \pm (\alpha _i ,\alpha _j )\hat{x}^\pm _{j,l}, \label{eq47}\\
[k_{i,0} ,x^\pm _{j,l} ] = \pm (\alpha _i ,\alpha _j )\hat{x}^\pm _{j,l}
, \qquad [k_{i,0} ,\hat{x}^\pm _{j,l} ] = \pm \widetilde{(\alpha _i ,\alpha _j)}x^\pm _{j,l},  \label{eq48}\\
k_{i,m + 1} = {1 \over 2}[h_{i + 1,1} , - h_{i - 1,1} ,k_{i,m}],\qquad  \tilde {h}_{i,1} = h_{i,1} + {1 \over 2}h_{i,0} ^2, \label{eq49}\\
\hat{x}^\pm _{i,2m + 1} = {1 \over 2}[h_{i + 1,1} , - h_{i -
1,1} ,\hat{x}^\pm _{i,2m} ],\qquad \hat{x}^\pm _{i,2m + 2} = {1
\over 2}[h_{i + 1,1} , - h_{i - 1,1} ,\hat{x}^\pm _{i,2m + 1} ], \label{eq50}\\
[\tilde {h}_{i,m + 1} ,x^\pm _{j,r} ] - [\tilde {h}_{i,m} ,x^\pm _{j,r + 1}] =
\pm {{(\alpha _{i,} \alpha _j )} \over 2}\{\tilde {h}_{i,m}
,x^\pm _{j,r} \} + (\pm \delta _{i,j + 1} - \delta _{i + 1,j} )\{k_{i,m}
,\hat{x}^\pm _{j,r} \}, \label{eq51}
\end{gather*}
we use before $\delta $ sign ``$+$'' in the case  $m + r \in 2Z_ + $, and sign
``$-$'' for $m + r \in 2Z_ + + 1$;
\begin{gather*}
[x^\pm _{i,m + 1} ,x^\pm _{j,r} ] - [x^\pm _{i,m} ,x^\pm _{j,r + 1} ] = \pm
{{(\alpha _{i,} \alpha _j )} \over 2}\{x^\pm _{i,m} ,x^\pm _{j,r}
\} + (\pm \delta _{i,j + 1} - \delta _{i + 1,j} )\{\hat{x}^\pm _{i,m}
,\tilde {x}^\pm _{j,r} \}, \nonumber
\end{gather*}
here signs before $\delta $ are defined also as in the previous formula;
\begin{gather*}
[\hat{x}^\pm _{i,m + 1} ,x^\pm _{j,r} ] - [\hat{x}^\pm _{i,m} ,x^\pm
_{j,r + 1} ] = \pm {\widetilde{(\alpha _i , \alpha _j)} \over
2}\{\hat{x}^\pm _{i,m} ,x^\pm _{j,r} \} \mp {\widetilde{(\alpha _i, \alpha _j)}
 \over 2}\{x^\pm _{i,m} ,\hat{x}^\pm _{j,r} \},\\ 
[\hat{x}^\pm _{i,m + 1} ,\hat{x}^\pm _{j,r} ] - [\hat{x}^\pm _{i,m}
, \hat{x}^\pm _{j, r + 1}] = \pm {{(\alpha _{i}, \alpha _j )}
\over 2}\{\hat{x}^\pm _{i,m} ,\hat{x}^\pm _{j,r} \} \mp
{\widetilde{(\alpha _{i,} \alpha _j)} \over 2}\{x^\pm _{i,m}, x^\pm _{j,r} \},\\ 
[k_{i,m + 1}, x^\pm _{j,r}] - [k_{i,m} ,x^\pm _{j, r+ 1} ] = \pm
{{(\alpha _{i,} \alpha _j )} \over 2}\{k_{i,m} ,x^\pm _{j,r} \}\\
\qquad{} +
(\pm \delta _{i,j + 1} - \delta _{i + 1,j} )\{(\bar {h}_{i,m} + \bar{h}_{i
+ 1,m} ),\hat{x}^\pm _{j,r} \},\\ 
[\tilde {h}_{i, 2m +1}, k_{j,r} ] = 2((\delta _{i,j} - \delta _{i,j + 1}
)\bar {k}_{i,2m + r + 1} + (\delta _{i,j} - \delta _{i, j-1}) \bar {k}_{i +
1, 2m + r + 1}), \\ 
[\tilde {h}_{i,2m} ,k_{j,2r + 1} ] = 0,\qquad
[k_{i,2k} ,k_{j,2l}] = 2(\delta _{i,j} - \delta _{i,j + 1}) \bar {h}_{i,2(k
+ l)} + 2(\delta _{i,j} - \delta _{i,j - 1} )\bar {h}_{i + 1, 2(k + l)}, \\ 
[k_{i,2m + 1} ,k_{j,2r} ] = 0, \qquad
\sum\limits_{\sigma \in S_3 } {[x^\pm _{i,\sigma (s_1 )}, [x^\pm _{i,\sigma
(s_2 )} ,x^\pm _{j, \sigma (s_3)} ]] = 0,} \\
\sum\limits_{\sigma \in S_3 } {[\hat{x}^\pm _{i,
\sigma (s_1 )} ,[x^\pm _{i,\sigma (s_2 )} ,x^\pm _{j, \sigma (s_3)} ]] =
0,} \qquad s_1 ,s_2 ,s_3 \in Z_ + . 
\end{gather*}
\end{theorem}

 Note that proof of this theorem is quite complicated and technical. We note only
 two issues of the proof:
 \begin{enumerate}\itemsep=0pt
 \item[1)] the formulas for even generators are proved similarly to the case of $Y(A(m,n))$ (see \cite{St});
 \item[2)] the def\/ining relations given in the Theorem~1 easily follows from those given in Theorem~2.
\end{enumerate}

In conclusion let us note that the problems of explicit description of
quantum double of Yangian of strange Lie superalgebra and
computation of the universal $R$-matrix are not discussed in paper and
will be considered in further work. This paper makes basis for such research.

\pdfbookmark[1]{References}{ref}
\LastPageEnding

\end{document}